\documentclass[graybox]{svmult}

\usepackage{type1cm}        
\usepackage{makeidx}         
\usepackage{graphicx}       
\usepackage{multicol}        
\usepackage[bottom]{footmisc}
\usepackage{newtxtext}       
\usepackage[varvw]{newtxmath}
\usepackage{url}
\makeindex 

\begin{document}

\title*{Renovating Calculus through Interdisciplinary Partnerships Using the SUMMIT-P Model
} 
\titlerunning{Interdisciplinary Collaborations: SUMMIT-P}

\author{Suzanne Dor\'ee and Jody Sorensen}

\institute{
Suzanne Dor\'ee \at Department of Mathematics, Statistics, and Computer Science\\ Augsburg University\\ Minneapolis, MN 55454\\
\email{doree@augsburg.edu}
\and Jody Sorensen \at Department of Mathematics, Statistics, and Computer Science\\ Augsburg University\\ Minneapolis, MN 55454\\
\email{sorensj1@augsburg.edu}
}

\maketitle

\abstract{This review paper highlights research findings from the authors’ participation in the SUMMIT-P project, which studied how to build and sustain multi-institutional interdisciplinary partnerships to design and implement curricular change in mathematics courses in the first two years of college, using the Curriculum Foundations Project (CFP) as a launchpad.  The CFP interviewed partner discipline faculty to learn about the mathematical needs of their students and how they use mathematics in their courses.  This paper summarizes research findings from the CFP and the SUMMIT-P project, and presents a detailed example of how these findings were implemented in the calculus sequence at Augsburg University to improve course focus, increase the relevance of course content, and provide opportunities for student to practice transference of the calculus to disciplinary contexts. 
This paper is based on the talk “Applied and Active Calculus Built Through Interdisciplinary Partnerships” presented at the 2022 AWM Research Symposium in the Session on “Research on the First Two Years of College Mathematics”.
} 

\section{Introduction}

In \emph{Leveraging Interdisciplinary Partnerships to Create an Impactful STEM Curriculum} \cite{GanterBourdeau2022}, Ganter and Bourdeau outline the critical role mathematics in the first two years of college plays in undergraduate STEM education and echo the call from multiple national reports over the past several decades to create a more relevant and engaging curriculum. They explain that 

\begin{quote}
     \ldots Interdisciplinary faculty partnerships are key in the development and implementation of a mathematics curriculum that allows students to really learn mathematical concepts while moving seamlessly to concepts in other disciplines. (p.1)
\end{quote}

This review paper highlights research findings from the Mathematical Association of America's (MAA) Curriculum Foundations project (CFP) that led to this call for interdisciplinary partnerships to work on introductory mathematics curriculum and the subsequent Synergistic Undergraduate Mathematics via Multi-institutional Interdisciplinary
Teaching Partnerships (SUMMIT-P) project that investigated implementation strategies.\footnote{This material is based upon work supported by the National Science Foundation under Grant No. 1822451. Any opinions, findings, and conclusions or recommendations expressed in this material are those of the author(s) and do not necessarily reflect the views of the National Science Foundation.}  

We explore these findings through a case study of a project at one of the SUMMIT-P participating institutions, Augsburg University, where we worked on renovating calculus through interdisciplinary partnerships. For an easy-to-share summary of the project see \cite{Doree2022}.

\section{The Curriculum Foundations Project, 1999-2011}

The MAA Committee on the Undergraduate Program in Mathematics (CUPM) subcommittee on Curriculum Renewal Across the First Two Years (CRAFTY) works ``To monitor ongoing developments with the intention of making general recommendations concerning the first two years of collegiate mathematics.'' \cite{CRAFTY} In 1999, members Susan Ganter and William Barker worked with CRAFTY to launch the Curriculum Foundations Project (CFP).  

With funding from participating institutions, the CFP held 22 workshops, each with a different partner discipline.  Each workshop brought together partner discipline faculty and mathematicians from many different institutions. Mathematics faculty listened to partner discipline faculty describe the mathematics that students in their disciplines needed to know. Subsequent conversations between mathematics and partner discipline faculty served to clarify terminology and create a shared understanding.  After each workshop, the partner discipline faculty summarized their recommendations in a report. The MAA published these reports and executive summaries in two volumes: The Curriculum Foundations Project: Voices of the Partner Disciplines \cite{GanterBarker2004a} and Partner Discipline Recommendations for Introductory College Mathematics and the Implications for College Algebra \cite{GanterHaver2011}.  

Although each discipline had their own ``wish list'' of mathematical skills and concepts, there was much agreement across the disciplines on meta-skills. The CFP findings \cite{GanterBarker2004b} showed that to prepare students for the partner discipline, mathematics courses should strive for depth over breadth in order to

\begin{quote}
\begin{itemize}
    \item Emphasize conceptual understanding \ldots, problem-solving skills [including transference] \ldots, mathematical modeling \ldots , [and] communication skills \ldots;
    \item Make the curriculum more appropriate for the needs of the partner disciplines \ldots;
    \item Encourage the use of active learning \ldots;
    \item Improve interdisciplinary cooperation \ldots; and
    \item Emphasize the use of appropriate technology.
    
    \hfill p.\ 3-7
\end{itemize}
\end{quote}

The next step was to implement this research in the curriculum at the national level.  The project team secured National Science Foundation funding (DUE-0511562) to study the implementation of the CFP findings in College Algebra. The results reported in \cite{GanterHaver2011} show a mixed success.  

\section{The SUMMIT-P Project, 2016-22}

In spite of mixed results in reforming College Algebra under the CFP, as well as previous challenges trying to reform Calculus nationally \cite{Bressoud2019}, Ganter et al. remained convinced that the CFP research had significant potential to influence curricular change nationally in mathematics in the first two years of college.

They turned again to the CRAFTY subcommittee in 2013 and proposed the SUMMIT-P project to study how to build and sustain interdisciplinary partnerships as a vehicle to improve the introductory mathematics curriculum.  These partnerships would be the key to overcoming the obstacles faced in both calculus and College Algebra reform efforts.

With funding from the National Science Foundation (lead grant DUE-1822451), the SUMMIT-P project (2016-2023) set out to 

\begin{quote}
\begin{itemize}
    \item Implement major recommendations from MAA’s Curriculum Foundations Project for the purpose of broadening participation in and institutional capacity for STEM learning, especially relative to teaching and learning in undergraduate mathematics courses;
     \item Foster a network of faculty and programs in order to promote community and institutional transformation through shared experiences and ideas for successfully creating functional interdisciplinary partnerships within and across institutions;
     \item Change the undergraduate mathematics curriculum in ways that support improved STEM learning for all students while building the STEM workforce of tomorrow; and
     \item Monitor various aspects of the Curriculum Foundations recommendations being implemented at participating institutions while measuring the impact on faculty and students. 
     
     \hfill \cite{Ganter2019a}, p.\ 770
\end{itemize} 
\end{quote}

Each of the ten SUMMIT-P participating institutions identified a mathematics course (or set of courses) in the first two years and committed to using the CFP findings and project-wide protocols to revise the course(s).  These revised courses would increase applied content and activities to better prepare students for partner discipline courses. 

Each institution built a team led by a Principal Investigator (PI) who was a mathematician with at least one Co-PI from a partner discipline.  Institutional teams ranged from three to five faculty members. A key element of the SUMMIT-P model was the creation of these interdisciplinary faculty learning communities (FLC). Filipas, et al. \cite{Filipas2022} outline the research basis for FLC as an ingredient in the Network Improvement Community theory of change. Bishop, et al. \cite{Bishop2020} describe the impact of using interdisciplinary FLC within the SUMMIT-P Project.  They explain
\begin{quote}
Academics often teach and research in silos. Moving beyond those bounds and interacting with so many different people from different disciplines was both instructive and joyful. While developing a shared understanding of course content across disciplines was no easy feat, it proved to be one of the most beneficial outcomes of the FLC. (p.\ 81)
\end{quote}

The SUMMIT-P project developed a variety of protocols that were used at each participating institution to build and sustain interdisciplinary teams, including
\begin{itemize}
\item Fishbowl discussions: mathematicians listen while partner discipline faculty discuss the CFP recommendations relevant to their discipline;  \cite{Filipas2022}, \cite{Hofrenning2020}
\item Partner discipline “wish lists” of content;
\item Follow-up conversations with partner discipline faculty; and
\item Review of partner discipline textbooks.
\end{itemize}
From this internal research, each institution revised their courses, including creating applied activities that connected the mathematics to the partner discipline. Bowers, et al.\ \cite{Bowers2020} describe how this process using a ``Plan-Do-Study'' cycle worked at several SUMMIT-P institutions.  Institutions also adapted activities developed by other institutions in the project. \cite{May2020}

Additionally, the SUMMIT-P project developed cross-institutional protocols including annual in-person meetings, frequent virtual meetings, professional development activities, and collaborative dissemination. A highlight of the collaborative work was that each SUMMIT-P institution hosted two highly-structured site visits for a visiting SUMMIT-P team, project evaluators, and stakeholders on the host's campus \cite{Filipas2022}, \cite{Piercey2020}.

The core research of the SUMMIT-P project was to study the faculty teams.  The research team identified a number of factors that influenced each institution's ability to build and sustain interdisciplinary partnerships and to implement revised courses. Four levels of factors identified by Slate Young, et al.\ \cite{SlateYoung2022} include: student, faculty, department, and institution.  For example, critical faculty-level factors included faculty enthusiasm, peer support, and synergy among participants. Department-level factors included the departmental culture of collaboration as well as the level of administrative support in the department. In addition to institutional support and institutional culture of collaboration, the size of the institution was another institutional-level factor.

\section{Augsburg University}

Augsburg University is a small, private college in Minneapolis, Minnesota with a long tradition of serving a diverse student body.  Augsburg students include new immigrants living in or near the campus's Cedar-Riverside neighborhood, first-generation students, adult learners through Augsburg's Adult Undergraduate program, students with disabilities through Augsburg's CLASS program, and students in recovery from chemical dependency through Augsburg's STEP-UP program.  Augsburg is nationally recognized as LGBTQIA-friendly\footnote{Lesbian, Gay, Bisexual, Transgender, Queer, Intersex, and Asexual}.   Over the past decade the racial diversity of Augsburg's student body has increased dramatically, especially the number of BIPOC\footnote{Black, Indigenous, and People of Color} students.  For example, 58\% of day students identified as people of color in the 2022-23 academic year compared to around 25\% in the 2009-10 academic year. 

Many of the key factors for success identified by the SUMMIT-P research in \cite{SlateYoung2022} are present at Augsburg.  At the institutional level, Augsburg faculty members regularly engage in curricular revision as part of our teaching, demonstrating faculty enthusiasm and peer support. Efforts resulting in presentations and publications are acknowledged as Scholarship of Teaching and Learning (SoTL) within scholarship expectations for faculty tenure and promotion, which shows institutional support.  The Mathematics, Statistics, and Computer Science (MSCS) Department is multidisciplinary and we often work with faculty members in the sciences, economics, and business on research and curricula, as part of a culture of collaboration.  Being small and nimble has allowed us to be ``incubators'' for curricular and pedagogical change.  For example, we teach a contextualized developmental algebra course and many courses incorporate active and inquiry-based learning. The ethos among our department faculty is highly collaborative and peer support is strong. We value applied mathematics as is reflected in our curriculum and our faculty's expertise -- two of our five mathematicians are in applied areas.  Augsburg students often need extra help in order to succeed and faculty members within the department and across the university recognize and reward efforts to build courses that engage students and support their success.

Augsburg offers one version of calculus consisting of Calculus I, Calculus II, and Multivariable Calculus. Students taking calculus are often working towards a major in biochemistry, biology, chemistry, computer science, data science, economics, education (secondary licensure in mathematics), mathematics, or physics.  
The calculus courses meet meet Monday, Wednesday, and Friday (MWF) for 70 minutes per class. Calculus I and II have have an additional 100 minute lab period on Tuesday or Thursday. 

Before the SUMMIT-P project, the MWF class sessions were primarily taught using traditional lectures punctuated by occasional small group activities.  Most applications were introduced in the lab period, which also included lessons on the use of computational technology, such as Microsoft Excel or Mathematica. Over several decades, we had incrementally reduced the theoretical content of these courses and increased the relevance of the course material.  For example we had replaced some integration techniques with an introduction to differential equations.  As Augsburg's student population changed, we recognized the need to accelerate these efforts. The SUMMIT-P project offered us support and a mechanism to tackle this work.

Augsburg University was one of the original participating institutions in the SUMMIT-P project with funding from the National Science Foundation (DUE-1625142).  The interdisciplinary team at Augsburg consisted of three mathematicians: Suzanne Dor\'ee (PI), Jody Sorensen, and Pavel B\v el\'ik, along with chemist Joan Kunz and economist Stella Hofrenning. 

In recognition that our calculus courses were already successful (as measured by pass rates of 75-85\% and passing grades in subsequent courses), 
the Augsburg team framed our curricular revision as a ``renovation'' of the courses.  This vision proved useful in building trust with mathematics faculty members who were not on the team because it honored their previous work. 

It also helped bring together mathematics faculty that we worked together to build shared departmental goals for calculus -- that students would
\begin{itemize}
    \item Appreciate that calculus is useful -- \textbf{relevance};
    \item Be able to recognize and apply calculus in new contexts -- \textbf{transference};
    \item Have procedural fluency and conceptual understanding as needed\footnote{For further discussion of these terms see \cite{NRC2001}.}; and
    \item Want to take more mathematics.
\end{itemize} 

\section{Focusing on Core Calculus Concepts}

The Augsburg team immediately recognized the challenge of increasing relevance and opportunities to practice transference in the calculus courses because the courses already contained arguably too many skills and concepts.  Over the previous decade we had worked to gradually streamline some of the course content, but we found it difficult to imagine how to reduce the course content further. 

We began with facilitating four fishbowl-style conversations with partner discipline faculty from (1) biology and environmental science; (2) computer science, business and economics; (3) chemistry; and (4) mathematics and physics. All disciplines recommended meta-cognitive goals aligned with CFP findings listed earlier.  Following these conversations we collected ``wish lists'' from chemistry, economics, and physics faculty since they had identified specific content recommendations.  We also reviewed introductory and advanced textbooks from these disciplines.

This work helped us to eliminate some topics from our calculus courses because that skill or concept was no longer used by partner disciplines or subsequent mathematics courses.  For example, we discovered that no courses used partial fractions or trigonometric substitution for integration, so we stopped teaching those techniques. 

While partner discipline faculty identified most calculus concepts and skills as necessary, when we dug a deeper we learned that the partner disciplines often relied on quite basic concepts and examples. As a result, we were able to simplify course content.  For example, partner disciplines made regular use of the sine function and sometimes used the cosine and tangent functions, but they rarely used the cosecant, secant, or cotangent functions so we focused on the first three functions.  Similarly the partner disciplines often used basic functions such as quadratics or exponentials and simple combinations such as $f(x) = e^{0.2x}$, $h(x) = xsin(x)$ or $c(x) = x + \frac{1}{x}$.  We were able to reduce emphasis on artificially-complicated combinations of functions.

Sometimes these discussions revealed content that could be moved from one course to another. When we learned that the only partner discipline at Augsburg that used algebraic methods to find limits was physics, and knowing that any mathematics or physics major would complete all three semesters of calculus, we moved the algebraic treatment of limits from Calculus I to Calculus II. 

By eliminating, reducing, and moving some topics, we created room to incorporate topics that were more important to partner discipline faculty.  We now include a few days on differential equations and a day introducing partial derivatives in Calculus I, more days on differential equations in Calculus II, and extended time on vector methods in Multivariable Calculus.  

An unexpected side effect of our work happened when we learned that both chemistry and economics students needed the content of Multivariable Calculus but often stopped after two semesters (Calculus I-II). As a result, we changed the prerequisite of Multivariable Calculus to allow students to take it immediately after Calculus I.  Now students taking two semesters may take Calculus I-Multivariable Calculus instead. We had to rearrange only a few days of content to support this change -- introducing integration by parts in Calculus I and reviewing it in Multivariable Calculus and moving polar coordinates to Multivariable Calculus.

Since Augsburg's Calculus I course serves a broad range of majors, and we decided to center Calculus I on the big ideas that we felt a student should encounter even if they only took one semester of calculus:
\begin{itemize}
    \item Rate of change, including some differentiation rules;
    \item The concept of a limit;
    \item Calculus as a tool for approximation;
    \item Calculus as a tool for optimization; and
    \item The Fundamental Theorem of Calculus.
\end{itemize}

Our renovated Calculus II course studies limits and convergence of functions, sequences, infinite series, and power series; extends modeling skills to include differential equations, integrals, and power series; presents a few examples of solving differential equations and integration techniques; and revisits the concept of approximation to include using power series.

\section{Making Calculus Relevant}

Not surprisingly, making mathematical content relevant to students appears to be correlated with student interest, motivation, engagement, and achievement.  \cite{Harahiewicz2016}, \cite{Hulleman2016}, \cite{Johansen2023}, \cite{Koskinen2022}, \cite{Levya2022}
We want students to see that calculus is useful and relevant to their lives and future studies. Focusing the calculus curriculum gave us the room to increase the relevance of our courses.  We now include applied examples nearly every day. See Table 1 for examples of applications and mathematical topics. 

Rather than first teaching the theory and then seeing applications, we flipped the order and use the applications to motivate the day's theoretical material.  We now dedicate the first 15-20 minutes of every MWF class in Calculus I and II to small group work on an applied Exploratory Activity.  

\begin{table}[!t]
\caption{Activity Topics and Content}
\label{tab:1}       
\begin{tabular}{p{2cm}p{5.2cm}p{4cm}}
\hline\noalign{\smallskip}
Subject & Context & Calculus Topic   \\
\noalign{\smallskip}\svhline\noalign{\smallskip}
Biology & Indeterminate Growth in Fish & Intro to Diff Eqs \\
& SIR Disease Model & Concavity \\
& Lynx-Hare Model & Systems of Diff Eqs \\
Business & Bank Account with regular deposits &  Separation of Variables \\
Chemistry & Molecular Attraction (Morse Potential) & Instantaneous Rate of Change \\ 
Economics & Marginal Cost vs. Average Cost & Local Linearity \\
 & Supply and Demand & Optimization \\
 & Travel Expenses & Optimization \\
 & Technology Costs & Sequences \\
 Education/Psych & Language Learning Rate & Euler's Method \\ 
 & Test Completion Time & Probability \\
 Engineering & Wind Power Generation & Optimization \\
 & Water Drainage & Improper Integrals \\
Medicine & Blood Alcohol Level & Rates of Change (Day 1)\\
& Tumor Growth under Medication & FTC\\
& Gompertz Tumor Growth & Euler's Method \\
Music & Overtones & Harmonic Series \\
Physics & Airplane Travel & Derivatives \\
& Pendulum Motion & Trig Derivatives \\
& Newton's Law of Cooling & Newton's Method \\
Weather & Wind Chill & Partial Derivatives \\

\noalign{\smallskip}\hline\noalign{\smallskip}
\end{tabular}
\end{table}

Since our calculus students are interested in many different majors, it was important that they see applications beyond physics and engineering.  We began by choosing the Briggs et al.\ text \cite{BriggsCochran} which provided a wide range of applied examples, particularly from the biological sciences.  We also wanted students to connect calculus to their daily lives, so we developed activities about topics such as the weather, driving, and drinking alcohol.  

A fertile source of applied examples connected to student's majors came from our partner discipline faculty.  Kunz is a chemist who also teaches environmental science.  Hofrenning is an economist who also teaches business.  With their help we were able to add examples from those disciplines.  Some examples were created by the Augsburg team, but many examples were adapted from textbooks with our partner discipline faculty member's help.  We also got ideas for examples from other institutions in the SUMMIT-P consortium and were able to share activities with faculty from other institutions.

While we were comfortable briefly introducing applied contexts, sometimes the examples from the partner discipline were too complicated for our broad audience.  For example, we considered creating an activity exploring the function that depicts how bond distance influences attraction, but the chemistry content was not sufficiently accessible.  We recast the example in terms of two people and their romantic interest in each other.  The example still introduces a function similar to what chemistry students will see later.

In addition to helping us find simple yet appropriate examples, partner discipline faculty helped us improve examples that we developed. As mathematicians, we are accustomed to seeing applied examples in our textbooks and often believe that is how the mathematics is used in applied contexts.  Partner discipline faculty were able to confirm or correct such instructor beliefs.  For example, in considering classic supply and demand curves, we mathematicians focused on finding the equilibrium price whereas an economist might focus on elasticity -- how a change in supply or demand affects the equilibrium price, which is an even more interesting mathematical question.

\section{Practicing Transference}

Since most students in our calculus courses use the skills and concepts they learn within courses in another major, it is vital that our calculus students develop the skill of transference. This desire was highlighted by our conversations with partner discipline faculty.  
According to the literature review by Perin \cite{Perin2011},  
\begin{quote}
    Contextualization is thought to promote transfer of learning and improve the retention of information. When information is learned in a context similar to that in which the skills will actually be needed, the application of learning to the new context may be more likely.
\end{quote}
Thus our work to increase relevance by connecting course material to the partner disciplines has the potential to contribute to transfer. 

We decided to ask students to practice transference at least once a week by dedicating the weekly Calculus I and Calculus II lab periods to building students' transference. During the lab periods, students work in small groups on extended applied activities. 
To avoid revealing what mathematics was needed for each lab activity, we proposed a significant change -- instead of asking students to practice recently-learned mathematical skills and concepts, each lab activity uses mathematics from at least two weeks earlier in the course or even an earlier course.  This delay has the added benefit of building in necessary review--  the first few lab activities in Calculus I focus on ideas from Precalculus and the first few lab activities in Calculus II focus on Calculus I material.

A typical lab activity begins with background information and data from some applied context. We ask students to describe the situation and estimate answers to questions using the numbers or graphs provided.  This first step helps establish the reasonableness of answers and begins sense-making.  
 
Next we provide an equation (or two) modeling the situation and again ask students to describe or estimate, now bringing in analytic and algebraic techniques.  We take care to not tell students what the necessary mathematics is.  For example, we might ask how some quantity is changing without mentioning ``the derivative'' and not hinting at what derivative rules to use.  We might ask when something is smallest or largest without mentioning ``optimization''.

Throughout the lab period, the instructor helps students see similarities between problems in the lab and earlier work and asks probing questions to help groups if they become stuck but we resist telling students what to do. As Billing explains in his comprehensive literature review and analysis, \cite{Billing2007}
\begin{quote}
Transfer is promoted if learners are shown how problems resemble each other, if they are expected to learn to do this themselves, if they are aware of how to apply skills in different contexts, if attention is directed to the underlying goal structure of comparable problems, if examples are varied and are accompanied by rules or principles (especially if discovered by the learners), and if learners’ self-explanations are stimulated. \hfill (p.\ 512)
\end{quote}
Fuchs et al. \cite{Fuchs2003} showed that directly teaching transfer skills and can improve ``far-transfer,'' such as to another course.
 
We built the lab activities from a variety of sources and our partner discipline faculty checked the authenticity of each lab.  Several labs were modified from the supplementary materials in our textbook.  We were also able to reuse some existing labs by editing the wording and moving the timing later.  See Table 2 for examples of lab activity applications and mathematical topics.  
\begin{table}[!t]
\caption{Lab Topics and Content}
\label{tab:2}      
\begin{tabular}{p{2cm}p{5.3cm}p{4cm}}
\hline\noalign{\smallskip}
Subject & Context & Calculus Topic   \\
\noalign{\smallskip}\svhline\noalign{\smallskip}
Biology & Foraging Animals & Optimization \\
Chemistry & Titration & Concavity \\
& Enzyme Kinetics & Quotient Rule \\
Design & Movie Theater Viewing Angle & Derivative Rules \\
& Time to Failure & Probability \\
Economics & Trend Adoption & Derivative Rules \\
& Elasticity & Definition of Derivative \\
& Gini Index of Wealth Distribution & Area \\
Earth Science & Climate Modeling & Linear Approximation \\
& Tornado Frequency & Exponential Functions \\
Geology & Richter Scale & Logarithms \\
Medicine & Periodic Drug Dosing & Geometric Series \\
Physics & Landing an Airplane & Chain Rule \\
& Designing a Water Clock & Volume of Revolution \\
& Position, Velocity, Acceleration, and Jerk & Taylor Polynomials \\

\noalign{\smallskip}\hline\noalign{\smallskip}
\end{tabular}
\end{table}

Two of our labs are published in an MAA Notes Volume: \cite{Belik2022} and \cite{Sorensen2022}. The first is based on using titration in chemistry to determine pH.  Students work with a data set and explore the connection between the second derivative and the ``equivalence point'' of the titration.  The second explores the idea of trend adoption from economics.  Again students explore a data set, this time showing the adoption rates of technology like microwaves and e-book readers.  The students fit a logistic model to a trend and use derivative rules to find rates of change.

\section{Augsburg Findings}

Because of the SUMMIT-P project and the work of our team at Augsburg University, we created a lively, applications-rich curriculum relevant to our diverse group of students that increased student engagement and lay the foundation for students to use their calculus knowledge in a wide range of majors and applied areas. We also built connections between mathematics faculty and faculty from partner disciplines across campus.  

In the five years that we have taught renovated calculus curriculum we have seen several encouraging outcomes.  Students no longer ask the the dreaded ``what is this good for?'' question, and they no longer complain about the time required for weekly lab periods.  In course surveys, students report enjoying the varied and relevant applications, say that they speak enthusiastically about the class to their peers, and recognize the importance of collaboration. Without direct prompting, over half of Calculus I students made comments about appreciating the utility of calculus on end-of-term course evaluations.  Students also expressed a strong sense of belonging and belief that they could succeed in mathematics courses.  

As we had hoped, instructors have observed that incorporating topics that students find important or relevant has heightened student interest and their engagement and motivation to learn the mathematics.  We have also found that using contexts that are familiar to students has helped students make sense of the mathematics and deepen their understanding. 

\section{Additional SUMMIT-P and MAA Resources}

Throughout the SUMMIT-P project, teams worked to embed a breadth of partner discipline co-developed applied activities.  For a sampling of twenty activities developed by SUMMIT-P institutions see \cite{Piercey2022}. Each article follows a template that provides information to the instructor, including background information on the partner discipline context, as well as the student version of the activity and solutions.  

Other institutions have described their journey through curricular revision in the SUMMIT-P Project including Embry Riddle Aeronautical University \cite{Wood2022}, Ferris State University \cite{Venkatesh2020}, LaGuardia Community College \cite{Lai2020}, Lee University \cite{Poole2022}, Norfolk State University \cite{Brucal-Hallare2020}, Oregon State University \cite{Beisiegel2020}, and Virginia Commonwealth University \cite{Ellwein2022}. 
Readers interested in revising curriculum in collaboration with the partner disciplines might find it useful to consult \cite{BeisiegelDoree2020} which includes a list of questions to consider in planning such a project as well as more information about the SUMMIT-P institutions. 

We also point readers to two key reports on calculus at the national level published by the MAA–– the MAA National Study of College Calculus \cite{Bressoud2015} and Addressing Challenges to the Precalculus to Calculus II Sequence through Case Studies \cite{Johnson2022} which also contain numerous references to the literature on calculus.  Each volume offers evidence-based suggestions for the improvement of calculus.

\bibliographystyle{unsrt}
\bibliography{DS_biblio_AWM}

\end{document}